\documentclass[10pt]{amsart}
\usepackage{amsfonts,amssymb,amscd,amsmath,enumerate,verbatim,calc}
\usepackage{mathrsfs}
\usepackage{xcolor}
\numberwithin{equation}{section}
\newtheorem{theorem}{Theorem}[section]
\newtheorem{corollary}[theorem]{Corollary}
\newtheorem{lemma}[theorem]{Lemma}
\newtheorem{proposition}[theorem]{Proposition}
\newtheorem{definition}[theorem]{Definition}
\newtheorem{
ples}[theorem]{Examples}

\DeclareMathOperator{\grad}{grad}
\DeclareMathOperator{\Hess}{Hess}
\DeclareMathOperator{\trace}{trace}
\DeclareMathOperator{\spec}{spec}

\begin{document}

\title[On the geometry of lift metrics and lift connections on the tangent bundle]
 {On the geometry of lift metrics and lift connections\\ on the tangent bundle}

\bibliographystyle{amsplain}
\author[ESMAEIL PEYGHAN, DAVOOD SEIFIPOUR, and ADARA M. BLAGA]{ESMAEIL PEYGHAN, DAVOOD SEIFIPOUR, and ADARA M. BLAGA}
\address{Department of Mathematics, Faculty of Science, Arak University,
Arak, 38156-8-8349, Iran.}
\email{e-peyghan@araku.ac.ir}
\address{Department of Mathematics, Abadan Branch, Islamic Azad University, Abadan, Iran.}
\email{d-seifipour@phd.araku.ac.ir}
\address{Department of Mathematics, Faculty of Mathematics and Computer Science, West University of Timi\c soara,
300223, Rom\^ ania.}
\email{adarablaga@yahoo.com
}



\begin{abstract}
We study lift metrics and lift connections on the tangent bundle $TM$ of a Riemannian manifold $(M,g)$. We also investigate the statistical and Codazzi couples of $TM$ and their consequences on the geometry of $M$. Finally, we prove a result on $1$-Stein and Osserman structures on $TM$, whenever $TM$ is equipped with the complete lift connection.
\end{abstract}

\maketitle

\section{Introduction}
The geometry of the tangent bundle with Riemannian lift metrics has been extensively studied in recent years (see \cite{ABG, AS, BPS, Davies, Dombrowski, GO, KS, O, YK1}, for instance). On the other hand, information geometry is an important and useful bridge between applicable and pure sciences, a combination between differential geometry and statistics \cite{Amari}. In this framework, methods of differential geometry are used and extended to probability theory. The mathematical point of view on information geometry was initiated by C. R. Rao. He showed that a statistical model could be a Riemannian manifold, via the Fisher information metric.
One of the main objects in this area are the statistical connections.
Statistical manifolds provide a geometrical model of probability distributions. The geometry of statistical manifolds has been applied to various fields of information science, information theory, neural networks, machine learning, image processing, statistical mechanics and etc. (\cite{Amari, AN, OPS, PFA, TVC}).
A statistical manifold is a differentiable manifold whose points are probability distributions (\cite{Amari, AN, Lauritzen, Matsuzoe}).
Precisely, a statistical structure on a differentiable manifold $M$ is a pair $(g,\nabla)$ such that $g$ is a Riemannian
metric and $\nabla$ is torsion-free affine connection with the property that $\nabla{g}$ is totally symmetric. A Riemannian
manifold $(M,g)$ together with Levi-Civita connection $\nabla$ of $g$ is a trivial example of statistical manifold. In other words, statistical manifolds can be regarded as generalizations of Riemannian manifolds.

In this paper we study the prolongations of statistical structures on manifolds to their tangent bundles with horizontal and complete lift connections.
We consider two Riemannian lift metrics on the tangent bundle $TM$ of a Riemannian manifold $(M,g)$, one of them is the twisted Sasaki metric $G^{f, h}$ (in particular, Sasaki metric) and the other one is the gradient Sasaki metric $g^f$.

The paper is organized as follows. After some preliminary considerations, in Section 3, we study the geometry of $TM$ equipped with the twisted Sasaki metric $G^{f,h}$ and the horizontal (respectively, complete) lift connection $\overset{H}{\nabla}$ (respectively, $\overset{C}{\nabla}$)
and we investigate some properties of the couples $(g^s, \nabla^{f,h})$ and $(g^{f_1}, \nabla^{f,h})$ on $TM$, where $\nabla^{f,h}$ is the Levi-Civita connection of the twisted Sasaki metric $G^{f,h}$, $g^s$ is the Sasaki metric and $g^{f_1}$ is the gradient Sasaki metric.
We also obtain some results on the lift to the tangent bundle of Killing vector fields and infinitesimal affine transformations.
In Section 4, we study the geometry of $TM$ equipped with the gradient Sasaki metric $g^f$ and the lift connection $\overset{C}{\nabla}$
and we investigate some properties of the couples $(g^s, \nabla^f)$ and $(G^{f,h}, \nabla^{f_1})$ on $TM$, where $\nabla^f$ is the Levi-Civita connection of the gradient Sasaki metric $g^f$.
We also study the necessary conditions for $(TM,g^s, \nabla^{f})$ and $(TM,G^{f,h}, \nabla^{f_1})$ to be Codazzi and statistical manifolds.
Finally, in Section 5, we prove a theorem on the spectral geometry of $TM$ and we deduce that $TM$ is globally Osserman, whenever it is equipped with the complete lift connection $\overset{C}{\nabla}$ and $\nabla$ is a flat connection.

\section{Preliminaries}

Let $\nabla$ be an affine connection on a differentiable manifold $M$, let $(x^i)$ be local coordinates on $M$ and let
$(x^i,y^i)$ be the induced coordinates on $TM$. Then $\{\frac{\partial}{\partial{x^{i}}}\mid_{(x,y)},\frac{\partial}{\partial{y^{i}}}\mid_{(x,y)}\}$ is the natural basis of $T_{(x,y)}TM$. It is known that, with respect to an affine connection, $T_{(x,y)}TM$ can be decomposed to $H_{(x,y)}TM\oplus V_{(x,y)}TM$, where $H_{(x,y)}TM$ is spanned by $\{\frac{\delta}{\delta{x^{i}}}\mid_{(x,y)}:=(\frac{\partial}{\partial x^i})^H=\frac{\partial}{\partial{x^{i}}}\mid_{(x,y)}-y^{k}\Gamma_{ki}^{j}(x)\frac{\partial}{\partial{y^{j}}}\mid_{(x,y)}\}$
and $V_{(x,y)}TM$ is spanned by $\{\frac{\partial}{\partial{y^{i}}}\mid_{(x,y)}:=(\frac{\partial}{\partial x^i})^V\}$, with $\Gamma_{ki}^{j}$ the connection coefficients of $\nabla$.
Denote by $\pi:TM\rightarrow M$, $\pi(x,y):=x$, and by $\chi(M)$ the set of all vector fields on $M$.

The various lifts of a vector field $X=X^{i}\frac{\partial}{\partial{x^{i}}}$ on $M$ (complete lift, horizontal lift and vertical lift, respectively) are defined as follows		 \begin{align*}
X^C=X^{i}\frac{\partial}{\partial{x^{i}}}+y^{a}\frac{\partial X^{i}}{\partial{x^{a}}}\frac{\partial}{\partial{y^{i}}},\ \ \
X^H=X^{i}\frac{\partial}{\partial{x^{i}}}-y^{a}\Gamma_{ai}^{k}X^{i}\frac{\partial}{\partial{y^{k}}},\ \ \
X^V=X^{i}\frac{\partial}{\partial{y^{i}}},
\end{align*}
(using Einstein summation convention).

According to \cite{YI}, the Lie brackets of the horizontal lift and vertical lift of vector fields are
\begin{align}\label{Lie bracket}
[X^H,Y^H]=[X,Y]^H-(R(X,Y)y)^V, \quad [X^H,Y^V]=(\nabla_{X}Y)^V-(T(X,Y))^V, \quad [X^V,Y^V]=0,
\end{align}
where $T$ is the torsion tensor field and $R$ is the curvature tensor field of $\nabla$.

The horizontal lift connection $\overset{H}{\nabla}$ and the complete lift connection $\overset{C}{\nabla}$ of the affine connection $\nabla$ are respectively defined by \cite{YI}:	 \begin{align}\label{horizontal conn}
\overset{H}{\nabla}_{X^{H}}Y^{H}=(\nabla_{X}Y)^{H}, \quad \overset{H}{\nabla}_{X^{H}}Y^{V}=(\nabla_{X}Y)^{V}, \quad \overset{H}{\nabla}_{X^{V}}Y^{H}=\overset{H}{\nabla}_{X^{V}}Y^{V}=0,\nonumber\\
\overset{C}{\nabla}_{X^{H}}Y^{H}=(\nabla_{X}Y)^{H}+(R(y,X)Y)^{V}, \quad
\overset{C}{\nabla}_{X^{V}}Y^{H}=\overset{C}{\nabla}_{X^{V}}Y^{V}=0,\\
\overset{C}{\nabla}_{X^{H}}Y^{V}=(\nabla_{X}Y)^{V},\quad
\overset{C}{\nabla}_{X^{C}}Y^{C}=(\nabla_{X}Y)^{C}, \quad \overset{C}{\nabla}_{X^{C}}Y^{V}=\overset{C}{\nabla}_{X^{V}}Y^{C}=(\nabla_{X}Y)^{V}\nonumber.
\end{align}	
It is known that $\nabla$ is flat and torsion-free if and only if $\overset{H}{\nabla}$ ($\overset{C}{\nabla}$) is torsion-free \cite{YI}.

For simplicity, in the rest of the paper, we shall write $\partial_{i}$, $\delta_{i}$ and $\partial_{\bar{i}}$ instead of $\frac{\partial}{\partial{x^{i}}}$, $\frac{\delta}{\delta{x^{i}}}$ and $\frac{\partial}{\partial{y^{i}}}$.

Let now $(M,g)$ be a Riemannian manifold. Similar to the lifts of vector fields and affine connections, we can define lifts of Riemannian metrics.

We construct the \textit{twisted Sasaki metric} $G^{f,h}$ on $TM$ as follows
\begin{align}\label{twisted}
G^{f,h}_{(x,y)}(X^{H},Y^{H})=f(x)g_{x}(X,Y), \quad
G^{f,h}_{(x,y)}(X^{V},Y^{H})=0, \quad  G^{f,h}_{(x,y)}(X^{V},Y^{V})=h(x)g_{x}(X,Y),
\end{align}
where $f ,h$ are strictly positive smooth functions on $M$. If $f=h=1$, then $G^{f,h}$ reduces to the Sasaki metric $g^s$ \cite{Sasaki}.

\begin{lemma}\cite{Belarbi-El Hendi 1}\label{twisted levi-Civita}
	Let $(M, g)$ be a Riemannian manifold, let $\nabla$ be the Levi-Civita connection of $g$ and let $(TM, G^{f,h})$ be its tangent bundle equipped with the twisted Sasaki metric. Then the Levi-Civita connection $\nabla^{f,h}$ of $G^{f,h}$ is given by
	\begin{align*}
	\nabla^{f,h}_{X^V}Y^{V}&=-\left(\frac{1}{2f}g(X,Y)\grad {h}\right)^H,\\
	\nabla^{f,h}_{X^V}Y^{H}&=\left(\frac{h}{2f}R(y,X)Y\right)^H+\left(\frac{Y(h)}{2h}X\right)^V,\\
	\nabla^{f,h}_{X^H}Y^{V}&=\left(\frac{h}{2f}R(y,Y)X\right)^H+\left(\frac{X(h)}{2h}Y+\nabla_{X}Y\right)^V,\\
	\nabla^{f,h}_{X^H}Y^{H}&=(\nabla_{X}Y+A_{f}(X,Y))^H-\frac{1}{2}(R(X,Y)y)^V,
	\end{align*}
	where $A_{f}(X,Y)=\frac{1}{2f}\Big(X(f)Y+Y(f)X-g(X, Y)\grad {f}\Big)$, $X, Y\in \chi(M)$ and $(x,y)\in TM$. In particular, the Levi-Civita connection $\nabla^s$ of the Sasaki metric $g^s$ is given by
	\begin{align*}
	\nabla^s_{X^H}Y^H&=(\nabla_{X}Y)^H-\frac{1}{2}(R(X,Y)y)^V, \quad \nabla^s_{X^V}Y^H=\frac{1}{2}(R(y,X)Y)^H,\\
	\nabla^s_{X^H}Y^V&=\nabla_{X}Y^V+\frac{1}{2}(R(y,Y)X)^H,
	\quad \nabla^s_{X^V}Y^V=0.
	\end{align*}
\end{lemma}

We construct also the \textit{gradient Sasaki metric} $g^f$ on $TM$ as follows
\begin{align}\label{gradient Sasaki metric}
g^f_{(x,y)}(X^H, Y^H)=g_{x}(X, Y), \quad g^f_{(x,y)}(X^H, Y^V)=0, \quad g^f_{(x,y)}(X^V, Y^V)=g_{x}(X, Y)+X_x(f)Y_{x}(f),
\end{align}
where $f$ is a strictly positive smooth function on $M$. If $f$ is a constant, then $g^f$ reduces to the Sasaki metric $g^s$.

\begin{lemma}\cite{Belarbi-El Hendi 2}\label{gradient levi-Civita}
Let $(M, g)$ be a Riemannian manifold, let $\nabla$ be the Levi-Civita connection of $g$ and let $(TM,g^f)$ be its tangent bundle equipped with the gradient Sasaki metric. Then the Levi-Civita connection $\nabla^f$ of $g^f$ is given by
	\begin{align*}
	\nabla^f_{X^V}Y^V&=-\frac{1}{2}X(f)(\nabla_{Y}\grad  f)^H-\frac{1}{2}Y(f)(\nabla_{X}\grad  f)^H,\\
	\nabla^f_{X^V}Y^H&=\frac{1}{2}(R(y,X)Y)^H+\frac{1}{2}X(f)(R(y,\grad  f)Y)^H+\frac{1}{2}X(f)(\nabla_{Y}\grad  f)^V\\
	&\quad +\frac{1}{2a}\{g(X, \nabla_{Y}\grad  f)-\frac{1}{2}Y(a)X(f)\}(\grad  f)^V,\\	
	\nabla^f_{X^H}Y^V&=\frac{1}{2}(R(y,Y)X)^H+\frac{1}{2}Y(f)(R(y, \grad  f)X)^H+\frac{1}{2}Y(f)(\nabla_{X}\grad  f)^V+(\nabla_{X}Y)^V\\
	&\quad +\frac{1}{2a}\{g(Y, \nabla_{X}\grad  f)-\frac{1}{2}X(a)Y(f)\}(\grad  f)^V,\\
	\nabla^f_{X^H}Y^H&=(\nabla_{X}Y)^H-\frac{1}{2}(R(X,Y)y)^V,	
	\end{align*}
where $a=1+\parallel \grad  f\parallel^2$, $X, Y\in\chi(M)$ and $(x,y)\in TM$.
\end{lemma}

\begin{definition}
Let $(M,g)$ be a Riemannian manifold and let $\nabla$ be an affine connection on $M$.
The pair $(g,\nabla)$ is said to be a \textit{Codazzi couple} on $M$ if the cubic tensor field $C:=\nabla{g}$, is totally symmetric, namely, the Codazzi equations hold:
\begin{align*}
(\nabla_{X}g)(Y,Z)=(\nabla_{Y}g)(Z,X)=(\nabla_{Z}g)(X,Y),
\end{align*}
for every $X, Y, Z \in \chi(M)$.
The triplet $(M,g,\nabla)$ is called a \textit{Codazzi manifold} and $\nabla$ is called a \textit{Codazzi connection}. Furthermore, if $\nabla$ is torsion-free, then $(M,g,\nabla)$ is a \textit{statistical manifold}, $(g,\nabla)$ is a \textit{statistical couple} and $\nabla$ is a \textit{statistical connection}.
\end{definition}

\section{Geometry of tangent bundle with twisted Sasaki metric}

In this section we study the geometry of $TM$ equipped with the twisted Sasaki metric (in particular the Sasaki metric).
\begin{definition}
Let $(M, g)$ be a Riemannian manifold and let $\nabla$ be an affine connection on $M$.
\begin{enumerate}
		\item[(1)] A vector field $X$ is said to be conformal (respectively, Killing) with respect to $g$, if $L_Xg=2\rho g$ (respectively, $L_{X}g=0$), where $\rho$ is a function on $M$ and the Lie derivative of $g$ in the direction of $X$ is given by
$(L_{X}g)(Y, Z):=Xg(Y,Z)-g(L_{X}Y,Z)-g(Y,L_XZ)$.
		\item[(2)] A vector field $X$ is said to be an infinitesimal affine transformation on $M$ with respect to $\nabla$, if $L_{X}\nabla=0$, where the Lie derivative of $\nabla$ in the direction of $X$ is given by $(L_{X}\nabla)(Y, Z):=L_{X}(\nabla_{Y}Z)-\nabla_{Y}(L_{X}Z)-\nabla_{[X,Y]}Z$.
\end{enumerate}
\end{definition}

Now we study conditions under which $X^V$ and $X^H$ are Killing vector fields for $G^{f,h}$.

By a direct computation and using (\ref{twisted}) and (\ref{Lie bracket}), we get
\begin{align*}
(L_{X^V}G^{f, h})(Y^V, Z^V)&=0,\\
(L_{X^V}G^{f, h})(Y^H, Z^V)
&=hg(\nabla_{Y}X-T(Y, X),Z),\\
(L_{X^V}G^{f, h})(Y^H, Z^H)&=0.
\end{align*}
If $\nabla$ is torsion-free, then $X^V$ is a Killing vector field for $G^{f,h}$  if and only if  $\nabla_{Y}X=0$. Moreover, using (\ref{twisted}) and (\ref{Lie bracket}), a straightforward computation gives
\begin{align*}
(L_{X^H}G^{f, h})(Y^V, Z^V)&=X(h)g(Y,Z)+h\Big((\nabla_{X}g)(Y,Z)+g(T(X,Y),Z)+g(Y,T(X,Z))\Big),\\
(L_{X^H}G^{f, h})(Y^H, Z^V)
&=hg(R(X,Y)y, Z),\\
(L_{X^H}G^{f, h})(Y^H, Z^H)
&=X(f)g(Y,Z)+f(L_{X}g)(Y,Z).
\end{align*}
If $\nabla$ is torsion-free, then $X^H$ is a Killing vector field for $G^{f,h}$  if and only if
\begin{align*}
(\nabla_{X}g)(Y,Z)=-\frac{X(h)}{h}g(Y, Z),\ \ R(X, Y)Z=0,\ \ (L_{X}g)(Y,Z)=-\frac{X(f)}{f}g(Y, Z), \ \ \ \forall \ Y, Z\in\chi(M).
\end{align*}
Thus we get the following

\begin{proposition}
Let $(M, g)$ be a Riemannian manifold and let $(TM, G^{f,h})$ be its tangent bundle equipped with the twisted Sasaki metric. Then the following assertions hold
\begin{enumerate}
	\item[(1)] if $\nabla$ is a torsion-free affine connection on $M$, then $X^V$ is a Killing vector field for $G^{f,h}$  if and only if  $X$ is a parallel vector field;
	\item[(2)] if $\nabla$ is a torsion-free affine connection on $M$, then $X^H$ is a Killing vector field for $G^{f,h}$  if and only if  $X$ is a conformal vector field on $(M, g)$ and $(\nabla_{X}g)(Y,Z)=-\frac{X(h)}{h}g(Y, Z)$, $R(X, Y)Z=0$, for all $Y, Z\in\chi(M)$;
	\item[(3)] if $\nabla$ is a torsion-free affine connection on $M$ and $h$ is constant, then $X^H$ is a Killing vector field for $G^{f,h}$  if and only if  $X$ is a conformal vector field on $(M, g)$, $\nabla$ is the Levi-Civita connection of $(M, g)$ and $R(X, Y)Z=0$, for all $Y, Z\in\chi(M)$;
	\item[(4)] if $\nabla$ is a torsion-free affine connection on $M$ and $f$ and $h$ are constant functions, then $X^H$ is a Killing vector field for $G^{f,h}$  if and only if  $X$ is a Killing vector field on $(M, g)$, $\nabla$ is the Levi-Civita connection of $(M, g)$ and $R(X, Y)Z=0$, for all $Y, Z\in\chi(M)$;
	\item[(5)] if $\nabla$ is the flat Levi-Civita connection on $(M, g)$ and $f$ and $h$ are constant functions, then $X^H$ is a Killing vector field for $G^{f,h}$  if and only if  $X$ is a Killing vector field on $(M, g)$.
\end{enumerate}	
\end{proposition}

Here, we compute the components of $\overset{H}{\nabla}G^{f,h}$ to study $(TM,G^{f,h},\overset{H}{\nabla})$. A direct computation gives
\begin{align}\label{L1}
(\overset{H}{\nabla}_{\delta_{i}}G^{f,h})(\delta_{j},\delta_{k})&=\delta_{i}G^{f,h}(\delta_{j},\delta_{k})-G^{f,h}(\overset{H}{\nabla}_{\delta_{i}}\delta_{j},\delta_{k})-G^{f,h}(\delta_{j},\overset{H}{\nabla}_{\delta_{i}}\delta_{k})\\
&=\delta_{i}(fg_{jk})-G^{f,h}((\nabla_{\partial_{i}}\partial_{j})^{H},(\partial_{k})^H)-G^{f,h}((\partial_{j})^{H},(\nabla_{\partial_{i}}\partial_{k})^{H})\nonumber\\
&=\partial_{i}(f)g_{jk}+f\partial_{i}(g_{jk})-fg(\nabla_{\partial_{i}}\partial_{j},\partial_{k})-fg(\partial_{j},\nabla_{\partial_{i}}\partial_{k})\nonumber\\
&=\partial_{i}(f)g_{jk}+f(\nabla_{\partial_{i}}g)(\partial_{j}, \partial_{k})\nonumber.
\end{align}
By a similar computation we get
\begin{align}\label{L2}
(\overset{H}{\nabla}_{\delta_{j}}G^{f,h})(\delta_{k},\delta_{i})=\partial_{j}(f)g_{ki}+f(\nabla_{\partial_{j}}g)(\partial_{k}, \partial_{i}),\ \ \  (\overset{H}{\nabla}_{\delta_{k}}G^{f,h})(\delta_{i},\delta_{j})=\partial_{k}(f)g_{ij}+f(\nabla_{\partial_{k}}g)(\partial_{i}, \partial_{j}).
\end{align}
We have also
\begin{align*}
(\overset{H}{\nabla}_{\partial_{\bar{i}}}G^{f,h})(\partial_{\bar{j}}, \partial_{\bar{k}})=0,\ \ \
(\overset{H}{\nabla}_{\delta_{i}}G^{f,h})(\delta_{j},\partial_{\bar{k}})=(\overset{H}{\nabla}_{\delta_{j}}G^{f,h})(\partial_{\bar{k}}, \delta_{i})=(\overset{H}{\nabla}_{\partial_{\bar{k}}}G^{f,h})(\delta_{i},\delta_{j})=0,
\end{align*}
\begin{align}\label{L3}
(\overset{H}{\nabla}_{\partial_{\bar{i}}}G^{f,h})(\partial_{\bar{j}},\delta_{k})=(\overset{H}{\nabla}_{\partial_{\bar{j}}}G^{f,h})(\delta_{k},\partial_{\bar{i}})=0, \quad (\overset{H}{\nabla}_{\delta_{k}}G^{f,h})(\partial_{\bar{i}}, \partial_{\bar{j}})=\partial_{k}(h)g_{ij}+h(\nabla_{\partial_{k}}g)(\partial_{i}, \partial_{j}).
\end{align}
If $(TM,G^{f,h},\overset{H}{\nabla})$ is a Codazzi manifold, then the second equation of (\ref{L3}) implies
$(\nabla_{\partial_{k}}g)(\partial_{i}, \partial_{j})=-\frac{1}{h}\partial_{k}(h)g_{ij}$. Setting this equation in (\ref{L1}) and using (\ref{L2}) we get
$(\partial_k(f)-\frac{f}{h}\partial_k(h))g_{ij}=0$, and consequently $h\partial_k(f)=f\partial_k(h)$. Thus we get the following

\begin{theorem}
   Let $(M, g)$ be a Riemannian manifold, let $\nabla$ be an affine connection on $M$ and let $(TM, G^{f,h})$ be its tangent bundle equipped with the twisted Sasaki metric. Then the following statements hold
	\begin{enumerate}
		\item [(1)] if $(TM,G^{f,h},\overset{H}{\nabla})$ is a Codazzi manifold, then $(\nabla_Zg)(X, Y)=-\frac{1}{f}Z( f)g(X, Y)$ and $f\grad h=h\grad f$, for all $X, Y, Z\in\chi(M)$. Moreover, $\overset{H}{\nabla}$ is compatible with $G^{f,h}$;
		\item [(2)] if $(TM,G^{f,h},\overset{H}{\nabla})$ is a statistical manifold, then $\nabla$ is flat, $(\nabla_Zg)(X, Y)=-\frac{1}{f}Z(Z)g(X, Y)$ and $f\grad h=h\grad f$, for all $X, Y, Z\in\chi(M)$. Moreover, $\overset{H}{\nabla}$ reduces to the Levi-Civita connection of $G^{f,h}$;
		\item [(3)] if $(TM,G^{f,h},\overset{H}{\nabla})$ is a statistical manifold and $h$ is a constant, then $f$ is constant, $\nabla$ is the Levi-Civita connection of $g$ and $\overset{H}{\nabla}$ reduces to the Levi-Civita of $G^{f,h}$;
		\item [(4)] if $\nabla$ is the Levi-Civita connection of $g$ and $f, h$ are constant functions, then $\overset{H}{\nabla}$ is compatible with $G^{f,h}$; in particular, if $\nabla$ is flat, then $\overset{H}{\nabla}$ reduces to the Levi-Civita connection of $G^{f,h}$.
	\end{enumerate}
\end{theorem}

Now we focus on $(TM, G^{f,h},\overset{C}{\nabla})$. A direct computation gives
\begin{align}
(\overset{C}{\nabla}_{\delta_{i}}G^{f,h})(\delta_{j},\delta_{k})&=
\delta_{i}G^{f,h}(\delta_{j},\delta_{k})-G^{f,h}(\overset{C}{\nabla}_{\delta_{i}}\delta_{j},\delta_{k})-G^{f,h}(\delta_{j},\overset{C}{\nabla}_{\delta_{i}}\delta_{k})\nonumber\\
&=\delta_{i}(fg_{jk})-G^{f,h}((\nabla_{\partial_{i}}\partial_{j})^{H}+(R(y,\partial_{i}),\partial_{j})^V,(\partial_{k})^H)\nonumber\\
&\quad -G^{f,h}((\partial_{j})^{H},(\nabla_{\partial_{i}}\partial_{k})^{H}+(R(y,\partial_{i})\partial_{k})^V)\nonumber\\
&=\partial_{i}(f)g_{jk}+f\partial_{i}(g_{jk})-fg(\nabla_{\partial_{i}}\partial_{j},\partial_{k})-fg(\partial_{j},\nabla_{\partial_{i}}\partial_{k})\nonumber\\
&=\partial_{i}(f)g_{jk}+f(\nabla_{\partial_{i}}g)(\partial_{j}, \partial_{k}).\nonumber
\end{align}
By a similar computation we get
\begin{align*}
(\overset{C}{\nabla}_{\delta_{j}}G^{f,h})(\delta_{k},\delta_{i})=\partial_{j}(f)g_{ki}+f(\nabla_{\partial_{j}}g)(\partial_{k}, \partial_{i}),\\ (\overset{C}{\nabla}_{\delta_{k}}G^{f,h})(\delta_{i},\delta_{j})=\partial_{k}(f)g_{ij}+f(\nabla_{\partial_{k}}g)(\partial_{i}, \partial_{j}).
\end{align*}
We have also
\begin{align*}
(\overset{C}{\nabla}_{\partial_{\bar{i}}}G^{f,h})(\partial_{\bar{j}}, \partial_{\bar{k}})=0,
\end{align*}
\begin{align}\label{3.10}
(\overset{C}{\nabla}_{\delta_i}G^{f,h})(\delta_j, \partial_{\bar{k}})=-hy^{s}R^{t}_{sij}g_{tk}, \quad
(\overset{C}{\nabla}_{\delta_j}G^{f,h})(\partial_{\bar{k}},\delta_i)=-hy^{s}R^{t}_{sji}g_{kt}, \quad (\overset{C}{\nabla}_{\partial_{\bar{k}}}G^{f,h})(\delta_i,\delta_j)=0,\\
(\overset{C}{\nabla}_{\partial_{\bar{i}}}G^{f,h})(\partial_{\bar{j}},\delta_{k})=(\overset{C}{\nabla}_{\partial_{\bar{j}}}G^{f,h})(\delta_{k},\partial_{\bar{i}})=0, \quad (\overset{C}{\nabla}_{\delta_{k}}G^{f,h})(\partial_{\bar{i}},\partial_{\bar{j}})=\partial_{k}(h)g_{ij}+h(\nabla_{\partial_{k}}g)(\partial_{i},\partial_{j}).\nonumber
\end{align}
If $(TM, G^{f,h}, \overset{C}{\nabla})$ is a Codazzi manifold, from (\ref{3.10}), we get $y^{s}R^{t}_{sji}=0$. Differentiating with respect to $y^t$, we obtain $R^{t}_{hji}=0$, so $\nabla$ is a flat connection. Thus we get the following

\begin{theorem}
   Let $(M, g)$ be a Riemannian manifold, let $\nabla$ be a torsion-free affine connection on $M$ and let $(TM, G^{f,h})$ be its tangent bundle equipped with the twisted Sasaki metric. Then the following statements hold
   \begin{enumerate}
	\item [(1)] if $(TM,G^{f,h},\overset{C}{\nabla})$ is a Codazzi (respectively, statistical) manifold, then $\nabla$ is flat, $(\nabla_Zg)(X, Y)=-\frac{1}{f}Z( f)g(X, Y)$ and $f\grad h=h\grad f$, for all $X, Y, Z\in\chi(M)$. Moreover, $\overset{C}{\nabla}$ is compatible with $G^{f,h}$ (respectively, $\overset{C}{\nabla}$ reduces to the Levi-Civita connection of $G^{f,h}$);
	\item [(2)] if $(TM,G^{f,h},\overset{C}{\nabla})$ is a statistical manifold and $h$ is a constant, then $f$ is constant, $\nabla$ is the Levi-Civita connection of $g$ and $\overset{C}{\nabla}$ reduces to the Levi-Civita of $G^{f,h}$;
	\item [(3)] if $\nabla$ is the Levi-Civita connection of $g$ and $f, h$ are constant functions, then $\overset{C}{\nabla}$ is compatible with $G^{f,h}$; in particular, if $\nabla$ is flat, then $\overset{C}{\nabla}$ reduces to the Levi-Civita connection of $G^{f,h}$.
\end{enumerate}
\end{theorem}

Here we provide conditions for the vertical vector field $X^V$ to be an infinitesimal affine transformation on $TM$ with respect to $\nabla^{f,h}$, whenever $(M,g)$ is a flat space and $\nabla$ is the Levi-Civita connection of $g$.

Using Lemma \ref{twisted levi-Civita} and (\ref{Lie bracket}), we get
\begin{align}
(L_{X^V}\nabla^{f,h})(Y^V, Z^V)&=\left(\frac{g(Y, Z)}{2f}\nabla_{\grad  h}X\right)^V,\label{D1}\\
(L_{X^V}\nabla^{f,h})(Y^H, Z^V)&=-\left(\frac{1}{2f}g(\nabla_{Y}X,Z)\circ\pi\right)(\grad  h)^H,\label{D2}\\
(L_{X^V}\nabla^{f,h})(Y^H, Z^H)&=-\left(\nabla_{A_{f}(Y,Z)+\nabla_{Y}Z}X+\frac{Y(h)}{2h}\nabla_{Z}X+\nabla_{Y}\nabla_{Z}X+\frac{Z(h)}{2h}\nabla_{Y}X\right)^V\label{D3}.
\end{align}
Let $X^V$ be an infinitesimal affine transformation on $TM$ with respect to $\nabla^{f, h}$. Then from (\ref{D2}) we get $\nabla_YX=0$, for all $Y\in\chi(M)$, or $\grad h=0$. In both cases, (\ref{D1}) vanishes. Also, in the first case, (\ref{D3}) vanishes and in the second case, (\ref{D3}) reduces to $\nabla_{A_{f}(Y,Z)+\nabla_{Y}Z}X+\nabla_{Y}\nabla_{Z}X=0$. Thus we have the following

\begin{proposition}
	Let $(M, g)$ be a flat Riemannian manifold, let $\nabla$ be the Levi-Civita connection of $g$, let $(TM, G^{f, h})$ be its tangent bundle equipped with the twisted Sasaki metric and let $\nabla^{f,h}$ be the Levi-Civita connection of $G^{f, h}$. Then $X^V$ is an infinitesimal affine transformation on $TM$ with respect to $\nabla^{f, h}$  if and only if  $X$ is parallel, or $h$ is constant and $\nabla_{A_{f}(Y,Z)+\nabla_{Y}Z}X+\nabla_{Y}\nabla_{Z}X=0$, for all $Y,Z\in\chi(M)$.
\end{proposition}

Here we provide necessary and sufficient conditions for the horizontal vector field $X^H$ to be an infinitesimal affine transformation on $TM$ with respect to $\nabla^{f,h}$.

Using Lemma \ref{twisted levi-Civita}, (\ref{Lie bracket}) and considering $Y^H(F^V)=(Y(F))^V$, $F^VY^V=(FY)^V$, for all $F\in C^\infty(M)$, we get
\begin{align*}
(L_{X^H}\nabla^{f,h})(Y^V,Z^V)&=\left(\frac{X(f)}{2f^2}g(Y,Z)\circ\pi\right)(\grad {h})^H-\left(\frac{1}{f}g(Y,Z)\circ\pi\right)[X, \grad {h}]^H\\
&=\frac{1}{f}\Big(\frac{X(f)}{2f}-1\Big)\Big(g(Y,Z)\circ\pi\Big)\Big(\grad h-[X, \grad {h}]\Big)^H,\\
(L_{X^H}\nabla^{f,h})(Y^H,Z^V)
&=(R(X,Y)Z)^V+X^H\left(\frac{Y(h)}{2h}\right)^VZ^V-\left(\frac{[X,Y](h)}{2h}\right)^VZ^V\\
&=\Big(R(X,Y)Z+X\left(\frac{Y(h)}{2h}\right)Z-\frac{[X,Y](h)}{2h}Z\Big)^V\\
&=\Big(R(X,Y)Z+\frac{1}{2h}\Big(-\frac{X(h)Y(h)}{h}+Y(X(h))\Big)Z\Big)^V,\\
(L_{X^H}\nabla^{f,h})(Y^H,Z^H)&=\Big((L_{X}\nabla)(Y,Z)+[X, A_{f}(Y,Z)]+A_{f}(Y, [Z, X])-A_{f}([X,Y], Z)]\Big)^H.
\end{align*}

\begin{proposition}
Let $(M, g)$ be a flat Riemannian manifold, let $\nabla$ be the Levi-Civita connection of $g$, let $(TM, G^{f, h})$ be its tangent bundle equipped with the twisted Sasaki metric and let $\nabla^{f,h}$ be the Levi-Civita connection of $G^{f, h}$. Then $X^H$ is an infinitesimal affine transformation with respect to $\nabla^{f,h}$  if and only if  the following equations hold
\begin{align*}
\Big(X(f)-2f\Big)\Big(\grad h-[X, \grad {h}]\Big)=0,\ \ \ \ \ \
R(X,Y)Z=\frac{1}{2h}\Big(\frac{X(h)Y(h)}{h}-Y(X(h))\Big),
\end{align*}
\begin{align*}
(L_{X}\nabla)(Y,Z)+[X, A_{f}(Y,Z)]+A_{f}(Y, [Z, X])-A_{f}([X,Y], Z)]=0,
\end{align*}
for all $Y, Z\in\chi(M)$.
\end{proposition}

\begin{corollary}
Let $(M, g)$ be a flat Riemannian manifold, let $\nabla$ be the Levi-Civita connection of $g$, let $(TM, G^{f, h})$ be its tangent bundle equipped with the twisted Sasaki metric and let $\nabla^{f,h}$ be the Levi-Civita connection of $G^{f, h}$.
\begin{enumerate}
	\item [(1)] If $f$ is constant, then $X^H$ is an infinitesimal affine transformation with respect to $\nabla^{f,h}$  if and only if  $X$ is an infinitesimal affine transformation with respect to $\nabla$ and
	\begin{align*}
	[X, \grad {h}]=\grad h,\ \ \ \ \ \
	R(X,Y)Z=\frac{1}{2h}\Big(\frac{X(h)Y(h)}{h}-Y(X(h))\Big),
	\end{align*}
	for all $Y, Z\in\chi(M)$.
	\item [(2)] If $h$ is constant, then $X^H$ is an infinitesimal affine transformation with respect to $\nabla^{f,h}$  if and only if
	\begin{align*}
	(L_{X}\nabla)(Y,Z)+[X, A_{f}(Y,Z)]+A_{f}(Y, [Z, X])-A_{f}([X,Y], Z)]=0,
	\end{align*}
	and $R(X, Y)Z=0$, for all $Y, Z\in\chi(M)$.
	\item [(3)] If $f$ and $h$ are constant, then $X^H$ is an infinitesimal affine transformation with respect to $\nabla^{f,h}$  if and only if  $X$ is an infinitesimal affine transformation with respect to $\nabla$ and $R(X, Y)Z=0$, for all $Y, Z\in\chi(M)$.
\end{enumerate}
\end{corollary}

Here we provide necessary and sufficient conditions for the horizontal vector field $X^H$ to be an infinitesimal affine transformation on $TM$ with respect to $\overset{H}{\nabla}$, where $\nabla$ is an affine connection on $M$.

By a direct computation and using (\ref{Lie bracket}) and (\ref{horizontal conn}), we get
\begin{align*}
(L_{X^H}\overset{H}{\nabla})(Y^V, Z^V)&=0,\\
(L_{X^H}\overset{H}{\nabla})(Y^H, Z^V)&=\Big(R(X, Y)Z-T(X, \nabla_{Y}Z)+\nabla_{Y}T(X, Z)\Big)^V,\\
(L_{X^H}\overset{H}{\nabla})(Y^H, Z^H)&=((L_{X}\nabla)(Y,Z))^H-(R(X, \nabla_{Y}Z)y)^V+(\nabla_{Y}R(X, Z)y)^V.
\end{align*}
Thus we get the following

\begin{proposition}
Let $\nabla$ be a flat torsion-free connection on $(M,g)$. Then $X^H$ is an infinitesimal affine transformation on $TM$ with respect to $\overset{H}{\nabla}$  if and only if  $X$ is an infinitesimal affine transformation on $M$ with respect to $\nabla$.
\end{proposition}

Here we consider $(TM, g^s,\nabla^{f,h})$ as a statistical manifold. Direct computations give $(\nabla^{f,h}_{\partial_{\bar{i}}}g^s)(\partial_{\bar{j}},\partial_{\bar{k}})=0$ and
\begin{align}\label{IM1}
(\nabla^{f,h}_{\delta_{i}}g^s)(\delta_{j},\delta_{k})=-g(A_{f}(\partial_{i},\partial_{j}),\partial_{k})-g(\partial_{j}, A_{f}(\partial_{i},\partial_{k}))=-\frac{\partial_i(f)}{f}g_{jk},
\end{align}
\begin{align}\label{3.13}
(\nabla^{f,h}_{\delta_{i}}g^s)(\delta_{j},\partial_{\bar{k}})=\frac{1}{2}y^rR_{ijrk}\left(1-\frac{h}{f}\right), \quad (\nabla^{f,h}_{\partial_{\bar{k}}}g^s)(\delta_{i},\delta_{j})=-\frac{h}{f}y^rR_{ijrk},
\end{align}
\begin{align}\label{IM2}
(\nabla^{f,h}_{\partial_{\bar{i}}}g^s)(\partial_{\bar{j}},\delta_{k})=\frac{\partial_{k}(h)}{2f}g_{ij}-\frac{\partial_{k}(h)}{2h}g_{ij}, \quad (\nabla^{f,h}_{\delta_{k}}g^s)(\partial_{\bar{i}},\partial_{\bar{j}})=-\frac{\partial_{k}(h)}{h}g_{ij}.
\end{align}
If $(TM,g^s, \nabla^{f,h})$ is a statistical manifold, from (\ref{IM1}) and (\ref{IM2}) we deduce that $f$ and $h$ are constant. From (\ref{3.13}), we get also $\frac{1}{2}y^rR_{ijrk}\left(1+\frac{h}{f}\right)=0$. Differentiating with respect to $y^t$, we obtain $-\frac{1}{2}R_{ijkt}\left(1+\frac{h}{f}\right)=0$.

According to the above case we get the following
\begin{theorem}
Let $(M, g)$ be a Riemannian manifold, $(TM, g^s)$ be its tangent bundle equipped with the Sasaki metric and let $\nabla^{f,h}$ be the Levi-Civita connection of twisted Sasaki metric $G^{f,h}$. If $(TM, g^s, \nabla^{f,h})$ is a statistical manifold, then $f$ and $h$ are constant. Moreover, we have that $\nabla$ is flat or $f=-h$.
\end{theorem}

Now we study the necessary conditions for $(TM, g^{f_{1}},\nabla^{f,h})$ to be a statistical manifold. First we recall the definition of Hessian. The Hessian of a function $f\in C^\infty(M)$ taken with respect to an affine
connection $\nabla$ is the covariant derivative of the 1-form $df$, i.e.,
\begin{align*}
\Hess ^{f}(X, Y):=(\nabla df)(X, Y)=XY(f)-(\nabla_XY)(f),\ \ \ \forall \ X, Y\in\chi(M).
\end{align*}
It is worth noting that $\Hess ^{f}$ is symmetric if and only if $\nabla$ is torsion-free.

Direct computations give
\begin{align}
(\nabla^{f,h}_{\delta_{i}}g^{f_1})(\delta_{j},\delta_{k})&=-\frac{\partial_i(f)}{f}g_{jk}, \quad (\nabla^{f,h}_{\partial_{\bar{i}}}g^{f_1})(\partial_{\bar{j}},\partial_{\bar{k}})=0,\nonumber \\
\label{3.15}
(\nabla^{f,h}_{\partial_{\bar{i}}}g^{f_1})(\partial_{\bar{j}},\delta_{k})&=\frac{\partial_{k}(h)}{2f}g_{ij}-\frac{\partial_{k}(h)}{2h}\Big(g_{ij}+\partial_{j}(f_1)\partial_{i}(f_1)\Big), \\
(\nabla^{f,h}_{\delta_{k}}g^{f_1})(\partial_{\bar{i}},\partial_{\bar{j}})&=-\frac{\partial_{k}(h)}{h}\Big(g_{ij}+\partial_{j}(f_1)\partial_{i}(f_1)\Big) \nonumber\\
&\quad +\Hess ^{f_1}(\partial_k, \partial_i)\partial_j(f_1)+\Hess ^{f_1}(\partial_k, \partial_j)\partial_i(f_1).\label{3.16}
\end{align}
If $(TM,g^{f_{1}}, \nabla^{f,h})$ is a statistical manifold, from (\ref{3.15}) and (\ref{3.16}), we get
\begin{align*}
\frac{\partial_{k}(h)}{2f}g_{ij}+\frac{\partial_{k}(h)}{2h}\{g_{ij}+\partial_{j}(f_1)\partial_{i}(f_1)\}&=\Hess ^{f_1}(\partial_k, \partial_i)\partial_j(f_1)+\Hess ^{f_1}(\partial_k, \partial_j)\partial_i(f_1).
\end{align*}
Also we have
\begin{align}\label{3.17}
(\nabla^{f,h}_{\delta_{i}}g^{f_1})(\delta_{j},\partial_{\bar{k}})=\frac{y^r}{2}R_{ijrk}+\frac{y^r}{2}R_{ijr}^{s}\partial_s(f_1)\partial_{k}(f_1)-\frac{h}{2f}y^rR_{rkij}, \quad (\nabla^{f,h}_{\partial_{\bar{k}}}g^s)(\delta_{i},\delta_{j})=-\frac{h}{f}y^rR_{rkij}.
\end{align}
Since $(TM,g^{f_{1}}, \nabla^{f,h})$ is statistical, from (\ref{3.17}), we get
\begin{align*}
\frac{y^r}{2}\{R_{ijrk}+R_{ijr}^{ s}\partial_s(f_1)\partial_{k}(f_1)+\frac{h}{f}R_{rkij}\}=0.
\end{align*}
Differentiating with respect to $y^t$, we obtain
\begin{align*}
\left(1+\frac{h}{f}\right)R_{ijtk}+R_{ijt}^{ s}\partial_s(f_1)\partial_{k}(f_1)=0.
\end{align*}
Thus we get the following

\begin{theorem}
Let $(M, g)$ be a Riemannian manifold, let $\nabla$ be an affine connection on $M$, let $(TM, g^{f_1})$ be its tangent bundle equipped with the gradient Sasaki metric and let $\nabla^{f,h}$ be the Levi-Civita connection of the twisted Sasaki metric $G^{f,h}$. If $(TM, g^{f_1}, \nabla^{f,h})$ is a statistical manifold, then
\begin{align*}
&\Hess ^{f_1}(Z, X)Y(f_1)+\Hess ^{f_1}(Z, Y)X(f_1)=\frac{1}{2(f+h)}Z(h)g(X, Y)+\frac{1}{2h}Z(h)Y(f_1)X( f_1),
\end{align*}
and
\begin{align*}
\left(1+\frac{h}{f}\right)R(X, Y)Z=(R(Y, X)Z)(f_1)\grad(f_1),
\end{align*}
for all $X, Y, Z\in\chi(M)$.
\end{theorem}

\section{Geometry of tangent bundle with gradient Sasaki metric}

In this section we study the geometry of $TM$ equipped with the gradient Sasaki metric $g^f$.

First we study the necessary and sufficient conditions for the vector fields $X^V$ and $X^H$ to be Killing for $g^f$.

By a direct computation and using (\ref{Lie bracket}) and  (\ref{gradient Sasaki metric}), we get
\begin{align*}
(L_{X^V}g^f)(Y^V,Z^V)&=0,\\
(L_{X^V}g^f)(Y^H,Z^V)&=g(\nabla_{Y}X,Z)+(\nabla_{Y}X)(f)Z(f)=g(\nabla_{Y}X,Z)+g((\nabla_{Y}X)(f)\grad(f), Z),\\
(L_{X^V}g^f)(Y^H,Z^H)&=0.
\end{align*}	

Using (\ref{Lie bracket}) and (\ref{gradient Sasaki metric}), straightforward computations give
\begin{align*}
(L_{X^H}g^f)(Y^V, Z^V)&=(\nabla_Xg)(Y,Z)+g(T(X,Y),Z)+g(Y, T(X, Z))\\
&\ \ \ +(T(X,Y))(f)Z(f)+(T(X,Z))(f)Y(f)\\
&\ \ \ +\Hess ^f(X, Y)Z(f)+\Hess ^f(X, Z)Y(f),\\
(L_{X^H}g^f)(Y^H, Z^V)&=g(R(X,Y)y, Z)+(R(X,Y)y)(f)Z(f)\\
&=g(R(X,Y)y, Z)+g((R(X,Y)y)(f)\grad(f), Z),\\
(L_{X^H}g^f)(Y^H, Z^H)&=(L_{X}g)(Y,Z).
\end{align*}	

Thus we get the following
\begin{proposition}
Let $(M, g)$ be a Riemannian manifold and let $(TM, g^f)$ be its tangent bundle equipped with the gradient Sasaki metric. Then the following assertions hold
\begin{enumerate}
\item[(1)] if $\nabla$ is a torsion-free affine connection on $M$, then $X^V$ is a Killing vector field for $g^f$,  if and only if  $\nabla
_YX=-(\nabla_YX)(f)\grad(f)$, for all $Y\in\chi(M)$;
\item[(2)] $X^H$ is a Killing vector field for $g^f$  if and only if  $X$ is a Killing vector field for $g$ and
	\begin{align*}
	(\nabla_Xg)(Y, Z)=-\Hess ^f(X, Y)Z(f)-\Hess ^f(X, Z)Y(f),\ \ \ \ R(X, Y)Z=-(R(X, Y)Z)(f)\grad(f),
	\end{align*}
\end{enumerate}
for all $Y, Z\in\chi(M)$.
  \end{proposition}

Here we compute the components of $\overset{C}{\nabla}g^{f}$ to study the  Codazzi and statistical structures for  $(TM, g^f,\overset{C}{\nabla})$.

A direct computation gives
\begin{align}\label{F1}
(\overset{C}{\nabla}_{\delta_{k}}g^f)(\partial_{\bar{i}}, \partial_{\bar{j}})&=\delta_{k}g^f(\partial_{\bar{i}}, \partial_{\bar{j}})-g^f(\overset{C}{\nabla}_{\delta_{k}}\partial_{\bar{i}}, \partial_{\bar{j}})-g^f(\partial_{\bar{i}}, \overset{C}{\nabla}_{\delta_{k}}\partial_{\bar{j}})\\
&=\delta_{k}g^f(\partial_{\bar{i}}, \partial_{\bar{j}})-g^f((\nabla_{\partial_{k}}\partial_{i})^V,\partial_{\bar{j}})-g^f(\partial_{\bar{i}},(\nabla_{\partial_{k}}\partial_{j})^V)\nonumber\\
&=\delta_{k}\{g_{ij}+\partial_{i}(f)\partial_{j}(f)\}-\{g(\nabla_{\partial_{k}}\partial_{i}, \partial_{j})+(\nabla_{\partial_{k}}\partial_{i})(f)\partial_{j}(f)\}\nonumber\\
&\quad -\{g(\partial_{i},\nabla_{\partial_{k}}\partial_{j})+\partial_{i}(f)(\nabla_{\partial_{k}}\partial_{j})(f)\}\nonumber\\
&=\partial_{k}(g_{ij})+\partial_{k}(\partial_{i}(f))\partial_{j}(f)+\partial_{i}(f)\partial_{k}(\partial_{j}(f))-g(\nabla_{\partial_{k}}\partial_{i}, \partial_{j})\nonumber\\
&\quad -(\nabla_{\partial_{k}}\partial_{i})(f)\partial_{j}(f)-g(\partial_{i},\nabla_{\partial_{k}}\partial_{j})-\partial_{i}(f)(\nabla_{\partial_{k}}\partial_{j})(f)\nonumber\\
&=(\nabla_{\partial_{k}}g)(\partial_{i},\partial_{j})+\partial_{k}(\partial_{i}(f))\partial_{j}(f)+\partial_{i}(f)\partial_{k}(\partial_{j}(f))\nonumber\\
&\quad -(\nabla_{\partial_{k}}\partial_{i})(f)\partial_{j}(f)-\partial_{i}(f)(\nabla_{\partial_{k}}\partial_{j})(f)\nonumber\\
&=(\nabla_{\partial_{k}}g)(\partial_{i},\partial_{j})+\partial_{j}(f)\Hess ^f(\partial_{k}, \partial_i)+\partial_{i}(f)\Hess ^f(\partial_{k}, \partial_j)\nonumber.
\end{align}
We have also
\begin{align}\label{F2}
(\overset{C}{\nabla}_{\partial_{\bar{i}}}g^f)(\partial_{\bar{j}}, \delta_{k})=(\overset{C}{\nabla}_{\partial_{\bar{j}}}g^f)(\delta_{k},\partial_{\bar{i}})=0,
\end{align}
\begin{align*}
(\overset{C}{\nabla}_{\delta_{i}}g^f)(\delta_{j},\partial_{\bar{k}})
&=-g(R(y,\partial_{i})\partial_{j},\partial_{k})-(R(y,\partial_{i})\partial_{j})(f)\partial_{k}(f),
\end{align*}
\begin{align*}
(\overset{C}{\nabla}_{\delta_{j}}g^f)(\partial_{\bar{k}},\delta_{i})=-g(R(y,\partial_{j})\partial_{i},\partial_{k})-(R(y,\partial_{j})\partial_{i})(f)\partial_{k}(f), \quad (\overset{C}{\nabla}_{\partial_{\bar{k}}}g^f)(\delta_{i},\delta_{j})=0,
\end{align*}
\begin{align}\label{4.1}
(\overset{C}{\nabla}_{\delta_{i}}g^f)(\delta_{j}, \delta_{k})=(\nabla_{\partial_{i}}g)(\partial_{j},\partial_{k}), \quad (\overset{C}{\nabla}_{\partial_{\bar{i}}}g^f)(\partial_{\bar{j}},\partial_{\bar{k}})=0.
\end{align}
Let $(TM, g^f, \overset{C}{\nabla})$ be a Codazzi manifold. The first equation of (\ref{4.1}) implies that $\nabla$ is Codazzi. On the other hand from (\ref{F1}) and (\ref{F2}) we deduce
\begin{align}\label{F4}
(\nabla_{\partial_{k}}g)(\partial_{i},\partial_{j})=-\partial_{j}(f)\Hess ^f(\partial_{k}, \partial_i)-\partial_{i}(f)\Hess ^f(\partial_{k}, \partial_j).
\end{align}
Since $\nabla$ is Codazzi, then the above equation gives
\begin{align}\label{F3}
\partial_{k}(f)\Hess ^f(\partial_{i}, \partial_j)-\partial_{i}(f)\Hess ^f(\partial_{k}, \partial_j)=T(\partial_i, \partial_k)(f)\partial_j(f).
\end{align}
Now, let $(TM, g^f, \overset{C}{\nabla})$ be a statistical manifold. It is known that $\overset{C}{\nabla}$ is torsion-free if and only if $\nabla$ is torsion-free and flat. So, we deduce that $\nabla$ is a statistical connection. In this case, (\ref{F3}) reduces to the following
\begin{align*}
\partial_{k}(f)\Hess ^f(\partial_{i}, \partial_j)-\partial_{i}(f)\Hess ^f(\partial_{k}, \partial_j)=0.
\end{align*}
Considering the above equation in (\ref{F4}) we get
\begin{align*}
(\nabla_{\partial_{k}}g)(\partial_{i},\partial_{j})=-2\partial_{k}(f)\Hess ^f(\partial_i, \partial_j).
\end{align*}
According to the above description we conclude the following

\begin{theorem}
Let $(M, g)$ be a Riemannian manifold, let $\nabla$ be an affine connection on $M$ and let $(TM, g^f)$ be its tangent bundle equipped with the gradient Sasaki metric. Then the following statements hold
\begin{enumerate}
	\item[(1)] if $(TM, g^f, \overset{C}{\nabla})$ is a Codazzi manifold, then $(M, g, \nabla)$ is a Codazzi  manifold,
	\begin{align*}
			(\nabla_Zg)(X, Y)=-\Hess ^f(Z, X)Y(f)-\Hess ^f(Z, Y)X(f),
	\end{align*}
	and
	\begin{align*}
			 R(X, Y)Z=-(R(X, Y)Z)(f)\grad(f),
	\end{align*}
such that	
\begin{align*}
		\Hess ^f(X, Y)Z(f)-\Hess ^f(Z, Y)X(f)=T(X, Z)(f)Y(f), \ \forall \ X, Y, Z\in\chi(M);
\end{align*}

	\item[(2)] if $(TM, g^f, \overset{C}{\nabla})$ is a statistical manifold, then $(M, g, \nabla)$ is a statistical manifold and
		\begin{align*}
				(\nabla_Zg)(X, Y)=-2\Hess ^f(X, Y)Z(f).
		\end{align*}
\end{enumerate}
\end{theorem}

Now we study the necessary conditions for $(TM, g^s,\nabla^{f})$ to be a statistical manifold. Direct computations give
\begin{align}
(\nabla^{f}_{\partial_{\bar{k}}}g^s)(\delta_{i},\delta_{j})&=-\frac{y^r}{2}R_{rkij}-\frac{y^r}{2}\partial_{k}(f)g(R(\partial_{r},\grad  f)\partial_{i},\partial_{j}) \nonumber\\
&\quad -\frac{y^r}{2}R_{rkji}-\frac{y^r}{2}\partial_{k}(f)g(\partial_{i},R(\partial_{r},\grad  f)\partial_{j}),\label{4.6}
\end{align}
\begin{align}\label{4.7}
(\nabla^{f}_{\delta_{i}}g^s)(\delta_{j},\partial_{\bar{k}})=-y^rR_{ijrk}-\frac{y^r}{2}\partial_{k}(f)g(\partial_{j},R(\partial_{r},\grad  f)\partial_{i}).
\end{align}
If $(TM,g^s, \nabla^f)$ is a statistical manifold, from (\ref{4.6}) and (\ref{4.7}), we get
\begin{align*}
 y^r\{R_{rkij}-\frac{\partial_{k}(f)}{2}g(\partial_{i},R(\partial_{r},\grad  f)\partial_{j})\}=0.
\end{align*}
Differentiating with respect to $y^t$, we obtain
\begin{align*}
R_{tkij}=\frac{\partial_{k}(f)}{2}g(\partial_{i},R(\partial_{t},\grad  f)\partial_{j}).
\end{align*}
We have also
\begin{align*}
(\nabla^{f}_{\delta_{i}}g^s)(\delta_{j},\delta_{k})= (\nabla^{f}_{\partial_{\bar{i}}}g^s)(\partial_{\bar{j}},\partial_{\bar{k}})=0,
\end{align*}
\begin{align}
(\nabla^{f}_{\partial_{\bar{i}}}g^s)(\partial_{\bar{j}},\delta_{k})&=\frac{\partial_{i}(f)}{2}g(\nabla_{\partial_{j}}\grad  f,\partial_{k})+\frac{\partial_{j}(f)}{2}g(\nabla_{\partial_{i}}\grad  f,\partial_{k})-\frac{\partial_{i}(f)}{2}g(\nabla_{\partial_{k}}\grad  f,\partial_{j}) \nonumber\\
&\quad -\frac{1}{2a}\{g(\partial_{i},\nabla_{\partial_{k}}\grad  f)-\frac{\partial_{k}(a)}{2}\partial_{i}(f)\}\partial_{j}(f),\label{4.10}
\end{align}
\begin{align}
(\nabla^{f}_{\delta_{k}}g^s)(\partial_{\bar{i}},\partial_{\bar{j}})&=-\frac{\partial_{i}(f)}{2}g(\nabla_{\partial_{k}}\grad  f,\partial_{j})-\frac{\partial_{j}(f)}{2}g(\nabla_{\partial_{k}}\grad  f,\partial_{i}) \nonumber\\
&\quad -\frac{1}{2a}\Big(g(\partial_{i},\nabla_{\partial_{k}}\grad  f)-\frac{\partial_{k}(a)}{2}\partial_{i}(f)\Big)\partial_{j}(f)-\frac{1}{2a}\Big(g(\partial_{j},\nabla_{\partial_{k}}\grad  f)\nonumber\\
&\quad -\frac{\partial_{k}(a)}{2}\partial_{j}(f)\Big)\partial_{i}(f).\label{4.11}
\end{align}
Since $(TM,g^s, \nabla^f)$ is statistical, from (\ref{4.10}) and (\ref{4.11}), we get
\begin{align*}
-\partial_{j}(f)g(\nabla_{\partial_{k}}\grad  f,\partial_{i})-\frac{1}{a}\Big(g(\partial_{j},\nabla_{\partial_{k}}\grad  f)-\frac{\partial_{k}(a)}{2}\partial_{j}(f)\Big)\partial_{i}(f)\\
=\partial_{i}(f)g(\nabla_{\partial_{j}}\grad  f,\partial_{k})+\partial_{j}(f)g(\nabla_{\partial_{i}}\grad  f,\partial_{k}).
\end{align*}

Thus we get the following
\begin{theorem}
	Let $(M, g)$ be a Riemannian manifold, $(TM, g^s)$ be its tangent bundle equipped with the Sasaki metric and let $\nabla^{f}$ be the Levi-Civita connection of the gradient Sasaki metric $g^f$. If $(TM, g^s, \nabla^{f})$ is a statistical manifold, then we have
	\begin{align*}
&-Y( f)g(\nabla_Z\grad f, X)-\frac{1}{a}\Big(g(\nabla_Z\grad f, Y)-\frac{1}{2}Z(a)Y(f)\Big)X( f)\\
&=X(f)g(\nabla_Y\grad f, Z)+Y(f)g(\nabla_X\grad f, Z),
	\end{align*}
	and
	\begin{align*}
R(X, Y)Z=\frac{1}{2}g(X, R(Z, \grad f)Y)\grad f,
	\end{align*}
	for all $X, Y, Z\in\chi(M)$.
\end{theorem}

Now we focus on $(TM, G^{f,h},\nabla^{f_1})$. Direct computations give
\begin{align}\label{E1}
(\nabla^{f_1}_{\delta_{i}}G^{f,h})(\delta_{j},\delta_{k})=\partial_{i}(f)g_{jk}, \quad (\nabla^{f_1}_{\partial_{\bar{i}}}G^{f,h})(\partial_{\bar{j}},\partial_{\bar{k}})=0,\quad
(\nabla^{f_1}_{\partial_{\bar{k}}}G^{f,h})(\delta_{i},\delta_{j})=0,
\end{align}
\begin{align}\label{4.2}
(\nabla^{f_1}_{\delta_{i}}G^{f,h})(\delta_{j},\partial_{\bar{k}})&=\frac{y^r}{2}\{(h-f)R_{ijrk}-f\partial_{k}(f_1)g(\partial_{j},R(\partial_{r},\grad  f_{1})\partial_{j})\}.
\end{align}
If $(TM,G^{f,h}, \nabla^{f_1})$ is a statistical manifold, from (\ref{E1}) we get $\partial_{i}(f)=0$, i.e., $f$ is constant.  (\ref{4.2}) implies
\begin{align*}
\frac{y^r}{2}\Big((h-f)R_{ijrk}-f\partial_{k}(f_1)g(\partial_{i},R(\partial_{r},\grad  f_{1})\partial_{j})\Big)=0.
\end{align*}
Differentiating with respect to $y^t$, we obtain
\begin{align*}
(h-f)R_{ijtk}-f\partial_{k}(f_1)g(\partial_{i},R(\partial_{t},\grad  f_{1})\partial_{j})=0.
\end{align*}
We have also
\begin{align}
(\nabla^{f_1}_{\partial_{\bar{i}}}G^{f,h})(\partial_{\bar{j}},\delta_{k})&=\frac{\partial_{i}(f_{1})}{2}fg(\nabla_{\partial_{j}}\grad  f_1,\partial_{k})+\frac{\partial_{j}(f_{1})}{2}fg(\nabla_{\partial_{i}}\grad  f_1,\partial_{k}) \nonumber\\
&\quad -\frac{\partial_{i}(f_{1})}{2}hg(\partial_{j},\nabla_{\partial_{k}}\grad  f_1)-\frac{h}{2a}\{g(\partial_{i},\nabla_{\partial_{k}}\grad  f_{1})-\frac{\partial_{k}(a)}{2}\partial_{i}(f_1)\}\partial_{j}(f_1),\label{4.4}
\end{align}
\begin{align}
(\nabla^{f_1}_{\delta_{k}}G^{f,h})(\partial_{\bar{i}},\partial_{\bar{j}})&=-\frac{\partial_{i}(f_1)}{2}hg(\nabla_{\partial_{k}}\grad  f_1,\partial_{j})-\frac{h}{2a}\{g(\partial_{i},\nabla_{\partial_{k}}\grad  f_{1})-\frac{\partial_{k}(a)}{2}\partial_{i}(f_1)\}\partial_{j}(f_1) \nonumber\\
&\quad -\frac{\partial_{j}(f_1)}{2}hg(\nabla_{\partial_{k}}\grad  f_1,\partial_{i})-\frac{h}{2a}\{g(\partial_{j},\nabla_{\partial_{k}}\grad  f_{1})-\frac{\partial_{k}(a)}{2}\partial_{j}(f_1)\}\partial_{i}(f_1)\nonumber\\
&\quad +\partial_{k}(h)g_{ij}.\label{4.5}
\end{align}
Since $(TM,G^{f,h}, \nabla^{f_1})$ is statistical, from (\ref{4.4}) and (\ref{4.5}), we get
\begin{align*}
&\frac{\partial_{i}(f_{1})}{2}fg(\nabla_{\partial_{j}}\grad f_1,\partial_{k})+\frac{\partial_{j}(f_{1})}{2}fg(\nabla_{\partial_{i}}\grad  f_1,\partial_{k})=\partial_{k}(h)g_{ij}
-\frac{h}{2a}\{g(\partial_{j},\nabla_{\partial_{k}}\grad  f_{1})\\
&-\frac{\partial_{k}(a)}{2}\partial_{j}(f_1)\}\partial_{i}(f_1)
-\frac{\partial_{j}(f_1)}{2}hg(\nabla_{\partial_{k}}\grad  f_1,\partial_{i}).
\end{align*}

Thus we get the following
\begin{theorem}
	Let $(M, g)$ be a Riemannian manifold, $(TM, G^{f,h})$ be its tangent bundle equipped with the twisted Sasaki metric and let $\nabla^{f_1}$ be the Levi-Civita connection of the gradient Sasaki metric $g^{f_1}$. If $(TM, G^{f,h},\nabla^{f_1})$ is a statistical manifold, then $f$ is constant. Moreover, we have
	\begin{align*}
	&fX(f_1)g(\nabla_Y\grad f_1, Z)+fY(f_1)g(\nabla_X\grad f_1, Z)=2g(X, Y)Z( h)\\
	&\ \ \ -hY(f_1)g(\nabla_Z\grad f_1, X)-\frac{h}{a}\Big(g(\nabla_Z\grad f_1, Y)-\frac{1}{2}Z(a)Y( f_1)\Big)X(f_1),
	\end{align*}
	and
	\begin{align*}
	(h-f) R(X, Y)Z=fg(X, R(Z, \grad f_1)Y)\grad f_1,
	\end{align*}
	for all $X, Y, Z\in\chi(M)$.
\end{theorem}

\section{$1$-Stein and Osserman structures on $TM$}

In this part we introduce two geometric concepts, such as $1$-Stein and Osserman space. Then we show that $TM$ is a $1$-Stein space whenever it is equipped with the complete lift connection $\overset{C}{\nabla}$. In the end we prove that if $\nabla$ is a flat connection, then $TM$ equipped with $\overset{C}{\nabla}$ is a globally Osserman space. It is known that $\overset{C}{\nabla}$ is the Levi-Civita connection of the complete lift metric on $TM$ (see \cite{YI} for more details). Studying $1$-Stein and Osserman structures on $TM$ with twisted Sasaki metric and gradient Sasaki metric are interesting ideas that can be studied in future.

The Jacobi operator $J_{X}(Y)=R(Y,X)X$ is a self-adjoint operator and it plays an important role in the curvature theory. Let $\spec \{J_{X}\}$ be the set of all eigenvalues of the Jacobi operator $J_{X}$ and $S(M,g)$ be the sphere bundle of unit tangent vector fields. One says that $(M,g)$ is \textit{Osserman at $p\in M$}, if for every $X, Y\in S_{p}(M,g)$, we have $\spec \{J_{X}\}=\spec \{J_{Y}\}$, i.e., the eigenvalues of $J_{X}$ are independent of the tangent vector at $p$. Furthermore, $(M,g)$ is \textit{pointwise Osserman}, if it is Osserman at each $p\in M$. Also, $(M,g)$ is \textit{globally Osserman} if, for any point $p\in M$ and any unit tangent vector $X\in T_{p}M$, the eigenvalues of the Jacobi operator depend neither on $X$ nor on $p$, i.e., the eigenvalues of $J_{X}$ are constant on $S(M,g)$. We recall that globally Osserman manifolds are clearly pointwise Osserman manifolds.

Let $(M,g)$ be a Riemannian manifold, $p \in M$, $Z \in S_{p}(M,g)$. Associated to the Jacobi operators, and natural number $t$, there exist some functions $f_{t}$ defined by $f_{t}(p,Z)=g(Z,Z)^{t}\trace (J_{Z}^{(t)})$, where $J_{Z}^{(t)}$ is the $t^{th}$ power of the Jacobi operator $J_{Z}$. We say that the Riemannian manifold $(M,g)$ is \textit{$k$-Stein at $p \in M$}, if $f_{t}(p,Z)$ is independent of $Z \in S_{p}(M,g)$ for every $1\leqslant t \leqslant k$. Moreover, $(M,g)$ is \textit{$k$-Stein} if it is $k$-Stein at each point.

\begin{lemma}\cite{GKVL}\label{remark1}
	Let $(M,g)$ be a $4$-dimensional Riemannian manifold. Then $(M,g)$ is pointwise Osserman if and only if $(M,g)$ is $2$-Stein.
\end{lemma}

Now we study the $1$-Stein and Osserman structure of $TM$, whenever it is equipped with the complete lift connection $\overset{C}{\nabla}$. We denote by $\bar{R}$ (respectively, $\bar{J}$) the Riemannian curvature tensor and the Jacobi operator of $\overset{C}{\nabla}$. Let $u\in TM$ and $v\in T_{u}(TM)$. So we have $v=X^k\delta_{k}+X^{\bar{k}}\partial_{\bar{k}}$, where $X^k, X^{\bar{k}}$ are smooth functions on $TM$. Direct computations give us
\begin{align*}
\bar{J}_{v}(\delta_{i})&=\bar{R}(\delta_{i},v)v=\bar{R}(\delta_{i},X^k\delta_{k}+X^{\bar{k}}\partial_{\bar{k}})(X^k\delta_{k}+X^{\bar{k}}\partial_{\bar{k}})\\
&=(X^{k})^2\bar{R}(\delta_{i}, \delta_{k})\delta_{k}+X^kX^{\bar{k}}\bar{R}(\delta_{i},\delta_{k})\partial_{\bar{k}}+X^{\bar{k}}X^k\bar{R}(\delta_{i},\partial_{\bar{k}})\delta_{k}+(X^{\bar{k}})^2\bar{R}(\delta_{i}, \partial_{\bar{k}})\partial_{\bar{k}}.
\end{align*}
Using (\ref{horizontal conn}) and a straightforward computation we get
\begin{align*}
\bar{R}(\delta_{i}, \delta_{k})\delta_{k}=(J_{\partial_{k}}(\partial_{i}))^H+\left(R(y,\partial_{i})\nabla_{\partial_{k}}\partial_{k}-R(y,\partial_{k})\nabla_{\partial_{i}}\partial_{k}+
\nabla_{\partial_{i}}R(y,\partial_{k})\partial_{i}-\nabla_{\partial_{k}}R(y,\partial_{i})\partial_{k}\right)^V.
\end{align*}
We have also
\begin{align*}
\bar{R}(\delta_{i}, \delta_{k})\partial_{\bar{k}}=(J_{\partial_{k}}(\partial_{i}))^V, \quad  \bar{R}(\delta_{i},\partial_{\bar{k}})\delta_{k}=\bar{R}(\delta_{i},\partial_{\bar{k}})\partial_{\bar{k}}=0.
\end{align*}
We set $A=R(y,\partial_{i})\nabla_{\partial_{k}}\partial_{k}-R(y,\partial_{k})\nabla_{\partial_{i}}\partial_{k}+
\nabla_{\partial_{i}}R(y,\partial_{k})\partial_{i}-\nabla_{\partial_{k}}R(y,\partial_{i})\partial_{k}$. Also from (\ref{horizontal conn}), we obtain $\bar{J}_{v}(\partial_{\bar{i}})=0$, thus the matrix representation of the Jacobi operator of $TM$ is
\begin{align*}
\bar{J}_{v}=\begin{bmatrix}
(X^k)^2(J_{\partial_{k}}(\partial_{i}))^H & 0\\
(X^k)^2A^V+(X^kX^{\bar{k}})(J_{\partial_{k}}(\partial_{i}))^V & 0\\
\end{bmatrix}
.
\end{align*}
Since $(X^k)^2(J_{\partial_{k}}(\partial_{k}))^H=0$, for $k=1,\dots, n$, we conclude that the principal diagonal entries are zero, so $\trace (\bar{J}_{v})=0$, i.e., $\trace$ is independent of $v$. Therefore $TM$ is a $1$-Stein space. Moreover, if $\nabla$ is a flat connection, then the Jacobi operator $\bar{J}_{v}$ of $TM$ equipped with $\overset{C}{\nabla}$ is zero. Thus $\spec \{\bar{J}_{V}\}=\{0\}$, so $TM$ is globally Osserman.
As mentioned above and using Lemma \ref{remark1}, we get the following

\begin{proposition}
	Let $(M, g)$ be a Riemannian manifold, let $\nabla$ be an affine connection and let $TM$ be its tangent bundle. Then the following statements hold
	\begin{enumerate}
		\item [(1)] $TM$ equipped with the complete lift connection $\overset{C}{\nabla}$ is a $1$-Stein space;
		\item [(2)] if $\nabla$ is a flat connection, then $TM$ equipped with $\overset{C}{\nabla}$ is globally Osserman. Moreover, if $M$ is $2$-dimensional, then $TM$ equipped with $\overset{C}{\nabla}$ is a $2$-Stein space.
	\end{enumerate}
\end{proposition}


\end{document}